\documentclass{amsart}
\newtheorem{Theorem}{Theorem}[section]
\newtheorem{Definition}[Theorem]{Definition}
\newtheorem{Proposition}[Theorem]{Proposition}

\newtheorem{Lemma}[Theorem]{Lemma}
\newtheorem{Corollary}[Theorem]{Corollary}
\theoremstyle{remark}

\newtheorem{Example}[Theorem]{Example}

\def\eps{\varepsilon}

\def\ovr{\overline}

\def\Om{\Omega}
\def\al{\alpha}
\def\gm{\gamma}

\def\Gm{\Gamma}
\def\Th{\Theta}

\def\dl{\delta}

\def\bd{\partial}
\def\lm{\lambda}

\def\si{\sigma}

\def\sm{\setminus}
\def\sbs{\subset}

\def\wtl{\widetilde}

\def\const{\operatorname{const}}

\def\sing{\operatorname{sing}}
\def\Aut{\operatorname{Aut}}

\def\reg{\operatorname{reg}}
\def\ind{\operatorname{ind}}

\def\re{{\mathbf {Re\,}}}

\def\be{\begin{enumerate}}
\def\ee{\end{enumerate}}
\def\bT{\begin{Theorem}}
\def\eT{\end{Theorem}}
\def\bP{\begin{Proposition}}
\def\eP{\end{Proposition}}
\def\bD{\begin{Definition}}
\def\eD{\end{Definition}}
\def\bE{\begin{Example}}
\def\eE{\end{Example}}
\def\bL{\begin{Lemma}}
\def\eL{\end{Lemma}}
\def\bC{\begin{Corollary}}
\def\eC{\end{Corollary}}
\def\A{{\mathcal A}}

\def\H{{\mathcal H}}

\def\oD{\ovr{\mathbb D}}
\def\aD{\mathbb D}

\def\aT{\mathbb T}
\def\aC{\mathbb C}

\def\aR{\mathbb R}
\def\aZ{\mathbb Z}

\def\C{{\mathcal C}}

\def\H{{\mathcal H}}

\def\R{{\mathcal R}}

\def\T{{\mathcal T}}

\begin{document}
\title{Fundamental group and analytic disks}
\author{Dayal Dharmasena and Evgeny A. Poletsky}
\begin{abstract} Let $W$ be a domain in a connected complex manifold $M$ and $w_0\in W$. Let $\A_{w_0}(W,M)$ be the space of all continuous mappings of a closed unit disk $\oD$ into $M$ that are holomorphic on the interior of $\oD$, $f(\bd\aD)\sbs W$ and $f(1)=w_0$. On the homotopic equivalence classes $\eta_1(W,M,w_0)$ of $\A_{w_0}(W,M)$ we introduce a binary operation $\star$ so that $\eta_1(W,M,w_0)$ becomes a semigroup and the natural mappings $\iota_1:\,\eta_1(W,M,w_0)\to\pi_1(W,w_0)$ and $\dl_1:\,\eta_1(W,M,w_0)\to\pi_2(M,W,w_0)$ are homomorphisms.
\par We show that if $W$ is a complement of an analytic variety in $M$ and if $S=\dl_1(\eta_1(W,M,w_0))$, then $S\cap S^{-1}=\{e\}$ and any element $a\in\pi_2(M,W,w_0)$ can be represented as $a=bc^{-1}=d^{-1}g$, where $b,c,d,g\in S$.
\par Let $\R_{w_0}(W,M)$ be the space of all continuous mappings of $\oD$ into $M$ such that $f(\bd\aD)\sbs W$ and $f(1)=w_0$. We describe its open dense subset $\R^{\pm}_{w_0}(W,M)$ such that any connected component of $\R^{\pm}_{w_0}(W,M)$ contains at most one connected component of $\A_{w_0}(W,M)$.
\end{abstract}
\thanks{The second author was partially supported by a grant from Simons Foundation.}
\keywords{holomorphic mappings, homotopic Oka principle, homotopy theory}
\subjclass[2010]{ Primary: 32Q55; secondary: 32H02, 32E30}
\address{Department of Mathematics, Faculty of Science, University of Colombo, Colombo 03, SRI LANKA} \email{dayaldh@sci.cmb.ac.lk}
\address{Department of Mathematics,  Syracuse University, \newline
215 Carnegie Hall, Syracuse, NY 13244}
\email{eapolets@syr.edu} 

\maketitle
\section{Introduction}
\par An {\it analytic disk} in a complex manifold $M$ is a continuous mapping $f$ of the closed unit disk $\ovr\aD$ into $M$ holomorphic on $\aD$.  We will denote the set of all such disks by $\A(M)$.  For a domain $W$ in $M$ we  introduce the space $\A(W,M)$ of all continuous mappings $f$ of the unit circle $\aT=\bd\aD$ into $W$ such that $f$ extends to a mapping $\hat f\in\A(M)$. If $w_0\in W$ then we denote by $\A_{w_0}(W,M)$ the subset of all $f\in\A(W,M)$ such that $f(1)=w_0$.
\par We let $\eta_1(W,M)$ to be the set of all connected components of $\A(W,M)$ and let $\eta_1(W,M,w_0)$ to be the set of all connected components of $\A_{w_0}(W,M)$. There is a natural mapping $\iota_1$ of the sets $\eta_1(W,M)$ or $\eta_1(W,M,w_0)$ into the sets $\pi_1(W)$ or $\pi_1(W,w_0)$ respectively.
\par In this paper we study the mapping $\iota_1$, its injectivity and its image. These questions originated in \cite{R1}, where L. Rudolph  showed that if $B$ is the braid group (the fundamental group of the complement of some set $A$ of hyperplanes in $\aC^n$), then $S=\iota_1((\eta_1(\aC^n\sm A,\aC^n,w_0)))$ is a semigroup, $S\cap S^{-1}=\{e\}$ and any element $a\in B$ can be represented as $a=bc^{-1}=d^{-1}g$, where $b,c,d,g\in S$. He called the elements of $S$ {\it quasipositive}.
\par In the recent paper \cite{KN} J. Koll\'ar and A. N\'emethi showed under additional assumptions that if $M$ is an algebraic variety with an isolated singularity $O$ and $W=M\sm O$, then $\iota_1:\,\eta_1(W,M)\to\pi_1(W)$ is an injection. They used this result to obtain more information about the singularity.
\par These results do not hold in general. For example, when $M={\mathbb{CP}^2}$ and $A$ is an algebraic variety in $M$ such that $\pi_1(M\sm A,w_0)=\aZ_p$ and $p$ is prime, then Rudolph's result evidently fails and the result of Koll\'ar and N\'emethi fails because the set $\eta_1(W,M)$ is infinite due to the homotopic invariance of the intersection index.
\par However Rudolph's result stays true if we change ingredients. We introduce on the set $\eta_1(W,M,w_0)$ a binary operation $\star$. With this operation $\eta_1(W,M,w_0)$ becomes a semigroup with unity. The natural mapping $\dl_1:\,\eta_1(W,M,w_0)\to\pi_2(M,W,w_0)$ is a homomorphism and we show that its image $S$ has the same properties as in Rudolph's result anytime when $W$ is a complement to an analytic variety in a connected complex manifold. When $\pi_1(M,w_0)=\pi_2(M,w_0)=0$ we obtain a complete analogy but in more general settings.
\par The problem of injectivity is more interesting and looks more difficult. To advance in this direction we consider the space $\R(W,M,w_0)$ of continuous mappings $f$ of $\oD$ into $M$ such that  $f(\aT)\sbs W$ and $f(1)=w_0$. We show that there is an open dense set $\R^\pm_{w_0}(W,M)$ in $\R_{w_0}(W,M)$ such that the natural mapping $\dl_1$ of $\eta_1(W,M,w_0)$ into the set $\rho^\pm_1(W,M,w_0)$ of all connected components of $\R^\pm_{w_0}(W,M)$ is an injection.
\par For future purposes we need to consider not domains $W\sbs M$ but Riemann domains $W$ over $M$. In Section \ref{S:bf} we prove basic facts about them. Since our constructions require more complicated compact sets than $\oD$, in Section \ref{S:ht} we introduce an operator $I_{K,\gm}$ that maps homotopic equivalence classes of holomorphic mappings of compact sets into homotopic equivalence classes of holomorphic mappings of the closed disk. The properties of this operator allows us in Section \ref{S:hfs} to introduce on $\eta_1(W,M,w_0)$ the structure of a semigroup. In Section \ref{S:phf} we establish major algebraic properties of $\eta_1(W,M,w_0)$. In particular, we obtain the description of the set $\eta_1(W,M)$ as the set of all $\pi_1$-conjugacy classes in $\eta_1(W,M,w_0)$.
\par In Section \ref{S:rho} we introduce the group $\rho_1(W,M,w_0)$ and prove its basic properties. In Section \ref{S:ctav} we consider the case when $W$ is the complement to an analytic variety in $M$ and $\Pi$ is an identity. The generalization of Rudolph's result is one of theorems in this section. The last Section \ref{S:cc} is devoted to the problem of injectivity.
\par We are grateful to Franc Forstneri\v c and Stephan Wehrli for valuable comments on the paper. We are also grateful to the referee whose corrections and commments significantly improved the exposition.
\section{Basic notions and facts}\label{S:bf}
\par In this paper $\aD(a,r)$ is an open disk of radius $r$ centered at $a$ and $\aT(0,r)$ is its boundary. We let $\aD=\aD(0,1)$ and $\aT=\aT(0,1)$.
\par A {\it Riemann domain} over a complex manifold $M$ is a pair $(W,\Pi)$, where $W$ is a connected Hausdorff complex manifold and $\Pi$ is a locally biholomorphic mapping of $W$ into $M$. Let $\hat d$ be a Riemann metric on $M$ and let $d$ be its lifting to $W$.
\par Let $K$ be a connected compact set in $\aC$ with connected complement. We denote by $\A(K,M)$ the set of all continuous mappings of $K$ into $M$ that are holomorphic on the interior $K^o$ of $K$. By $\A(K,W,M)$ we denote the set of all continuous mappings $f$ of $\bd K$ into $W$ such that there is a mapping $\hat f\in\A(K,M)$ coinciding with $\Pi\circ f$ on $\bd K$. The mapping $\hat f$ is unique. If $\zeta_0\in\bd K$ and $w_0\in W$ then the space $\A_{\zeta_0,w_0}(K,W,M)$ is the set of all $f\in\A(K,W,M)$ such that $f(\zeta_0)=w_0$.
\par The space $\T(W,M)$ consists of all pairs $(K,f)$, where $K\sbs\aC$ is a connected compact set with connected complement and $f\in\A(K,W,M)$. If $(K,f),(L,g)\in\T(W,M)$ we define the distance $d((K,f),(L,g))$ between  $(K,f)$ and $(L,g)$ as the sum of the Hausdorff distances  between the graphs of $f$ and $g$ on $\bd K$ and $\bd L$ respectively and between the graphs of $\hat f$ and $\hat g$ on $K$ and $L$ respectively. (The distance between points $(\zeta_1,w_1)$ and $(\zeta_2,w_2)$ on $\aC\times M$ is defined as $\max\{|\zeta_1-\zeta_2|,\hat d(w_1,w_2)\}$.) Since the graphs are compact, $d$ is a metric on $\T(W,M)$ and the topology on $\T(W,M)$ is induced by this metric. Clearly, this topology does not depend on the choice of $\hat d$.
\par The set $\T_{\zeta_0,w_0}(W,M)\sbs\T(W,M)$ consists of all pairs $(K,f)$ such that $f\in\A_{\zeta_0,w_0}(K,W,M)$. On this set and the sets $\A(K,W,M)$ and $\A_{\zeta_0,w_0}(K,W,M)$ we define the topology relative to the topology imposed on $\T(W,M)$.
\par Let $\T(M,M)$ be the set of pairs $(K,f)$, where $f\in\A(K,M)$. We define the mapping $\Pi_1$  of $\T(W,M)$ into the set $\T(M,M)$ as $\Pi_1(K,f)=(K,\hat f)$. Clearly, $\Pi_1$ is open and locally isometric.
\par Suppose that $K\sbs\aC$ is a compact set, $f\in\A(K,M)$ and the graph $\Gm_f$ of $f$ on $K$ has a Stein neighborhood $U$ in $\aC\times M$. Let $F$ be an imbedding of $U$ into ${\mathbb C}^N$ as a complex submanifold. By \cite[Theorem 8.C.8]{GR} there are
an open neighborhood $V$ of $F(\Gm_f)$ in ${\mathbb C}^N$ and a holomorphic retraction $P$ of $V$ onto $F(U)$.
\par Let $(L,g)$ be a pair, where $L\sbs \aC$ is a compact set, $g\in\A(L,M)$  and $\Gm_{g}\sbs U$. Then we let $\Phi(L,g)$ to be the pair $(L,h)$, where $h(\zeta)=F(\zeta,g(\zeta))$. Conversely, if $(L,h)$ is a pair, where $L\sbs \aC$ is a compact set, $h\in\A(L,\aC^N)$  and $\Gm_{h}\sbs V$, then we let $\Psi(L,h)$ to be the pair $(L,g)$, where $g=P_M\circ F^{-1}\circ P\circ h$ and $P_M$ is a projection of $\aC\times M$ onto $M$. Clearly, the mappings $\Phi$ and $\Psi$ are continuous and $\Psi\circ\Phi$ is the identity.
\par This construction leads to the following lemma.
\bL\label{L:hc} Let $K$ be a connected compact set in $\aC$ with connected complement. For every $\eps>0$ there is $\dl>0$ such that:\be
\item if $f\in\A(K,W,M)$ and pairs $(L,g_0)$ and $(L,g_1)$ lie in the $\dl$-neighborhood of $(K,f)$ in $\T(W,M)$, then there is a continuous path $(L,g_t)$ in the $\eps$-neighborhood of $(K,f)$ in $\T(W,M)$, $t\in[0,1]$, connecting $(L,g_0)$ and $(L,g_1)$. Moreover, if, additionally, a compact set $L'\sbs L$ and $g_0|_{L'}=g_1|_{L'}$, then we can assume that $g_t|_{L'}=g_0|_{L'}$ for all $t\in [0,1]$.
\item if $0\le t\le1$ and $(K,f_t)$, $(L_t,g_t)$, $w_t$ and $\xi_t\in L_t$ are continuous paths in $\T(W,M)$, $\T(W,M)$, $W$ and $\aC$ respectively and for all $0\le t\le 1$ the pairs $(L_t,g_t)$ lie in the $\dl$-neighborhood of $(K,f_t)$ and $d(g_t(\xi_t),w_t)<\dl$, then there is another continuous path $(L_t,h_t)$ in $\T(W,M)$ such that $h_t(\xi_t)=w_t$ and the pairs $(L_t,h_t)$ lie in the $\eps$-neighborhood of $(K,f_t)$ for all $0\le t\le1$. Moreover, if $g_t(\xi_t)=w_t$ for some $0\le t\le1$ then then we can assume that $h_t=g_t$.
\ee
\eL
\begin{proof} (1) It was shown in \cite[Theorem 3.1]{P} that the graph of $\hat f$ has a basis of Stein neighborhoods in $\aC\times M$. If $M=\aC^N$ then we connect $(L,\hat g_0)$ and $(L,\hat g_1)$ by the path $(L,\hat g_t)$, where
\[\hat g_t(\zeta)=(1-t)\hat g_0(\zeta)+t\hat g_1(\zeta).\]
\par In the general case, we take $\eps>0$ so that $\Pi_1^{-1}$ is defined on the $\eps$-neighborhood of $(K,\hat f)$. We choose $\dl>0$ so small that if we take pairs $\Phi(L,\hat g_0)$ and $\Phi(L,\hat g_1)$, connect them in $\aC^N$ by $(L,h_t)$ as above and  let $(L,\hat g_t)=\Psi(L,h_t)$, then the path $(L,\hat g_t)$ lies in the $\eps$-neighborhood of $(K,\hat f)$. Finally, we let $(L,g_t)= \Pi_1^{-1}(L,\hat g_t)$.
\par (2) For the proof of the second part we note that by \cite[Theorem 4.1]{P} the set $\wtl\Gm=\{(t,\zeta,\hat f_t(\zeta)):\,\zeta\in K, 0\le t\le1\}$ has a Stein neighborhood in $\aC\times\aC\times M$ and then the proof follows the same pattern as above.
\end{proof}
\par By part (1) of this theorem the spaces $\A(W,M)$ and $\A_{\zeta_0,w_0}(W,M)$ are locally path-connected and, therefore, their connected components are path-connected.
\section{Operator $I_{K,\gm}$}\label{S:ht}
\par Throughout this section $K$ will denote a connected  compact set in $\mathbb C$ with the connected complement. Let $\zeta_0\in\bd K$ and a base point $w_0\in W$. We say that $f,g\in\A_{\zeta_0,w_0}(K,W,M)$ are {\it $h$-homotopic} or $f\sim^h g$ if there is a continuous path connecting $f$ and $g$ in $\A_{\zeta_0,w_0}(K,W,M)$. The relation $\sim^h$ is evidently an equivalence and we denote the equivalence class of $f$ by $[f]_{\zeta_0,w_0}$ or $[f]$ if $\zeta_0$ and $w_0$ are fixed. The set of equivalence classes will be denoted by $\H_{\zeta_0,w_0}[K,W,M]$ or $\H_{\zeta_0,w_0}[K]$. It follows from Lemma \ref{L:hc}(1) that the equivalence classes are closed in $\A_{\zeta_0,w_0}(K,W,M)$.
\par Our goal is to construct a mapping of the set $\H_{\zeta_0,w_0}[K]$  into the set $\H_{1,w_0}[\oD]$. Firstly we do it when $K$ is the closure of a Jordan domain, i. e., $K$ is bounded by a Jordan curve (a homeomorphic image of a circle). Let $e_{1,\zeta_0}$ be a conformal mapping  of $\oD$ onto $K$ that maps 1 onto $\zeta_0$. We define the mapping $I_{K,\zeta_0}$ as $[f\circ e_{1,\zeta_0}]_{1,w_0}$. Since the group of conformal automorphisms of the unit disk with a fixed point on the boundary is connected, this mapping does not depend on the choice of $e$.
\par To define the mapping of $\H_{\zeta_0,w_0}[K]$ into $\H_{1,w_0}[\oD]$ for a general $K$ we will approximate $f\in\A_{\zeta_0,w_0}(K,W,M)$ by mappings on Jordan domains $\Om$ containing $K$. To determine a point in $\bd \Om$ where approximations are equal to $w_0$ we need an {\it access curve} to $K$ at $\zeta_0$, i.e., a continuous curve $\gm:\,[0,1]\to\aC$ such that $\gm(0)=\zeta_0$ and $\gm(t)\in\aC\sm K$ when $t>0$.
\par Let $\Om$ be a smooth Jordan domain containing $K$ whose boundary meets $\gm$. We let $\zeta_{\Om,\gm}=\gm(s_{\Om,\gm})$, where $s_{\Om,\gm}=\inf\{t:\,\gm(t)\in\bd \Om\}$. A pair $(\ovr \Om,g)\in\T(W,M)$ is an {\it $\eps$-approximation } of $(K,f)\in\T_{\zeta_0,w_0}(W,M)$ with respect to $\gm$ if $K\sbs \Om$, $g(\zeta_{\Om,\gm})=w_0$ and $(\ovr \Om,g)$ lies in the $\eps$-neighborhood of $(K,f)$.
\par The following proposition asserts the existence of $\eps$-approximations for every  $\eps>0$.
\bP\label{P:phiepsa} Let $f\in\A_{\zeta_0,w_0}(K,W,M)$ and let $\gm$ be an access curve to $K$ at $\zeta_0$. Then for every $\eps>0$ there is $\dl>0$ such that for any Jordan domain $\Om$ containing $K$ and lying in the $\dl$-neighborhood of $K$ and any point $\zeta\in\gm\cap \bd \Om$  there is a mapping $g\in\A_{\zeta,w_0}(\ovr \Om,W,M)$ such that the pair $(\ovr \Om,g)$ lies in the $\eps$-neighborhood of $(K,f)$.
\eP
\begin{proof}  Firstly, we prove this proposition when $f\in\A_{\zeta_0,w_0}(K,M,M)$ so $\hat f=f$. For the given $\eps>0$ we denote by $\eta$ the $\dl$ from Lemma \ref{L:hc}(2). The mapping $f$ is uniformly continuous on $K$. So there is $\dl_1>0$ such that $d( f(\zeta_1), f(\zeta_2))<\eta/2$ when $|\zeta_1-\zeta_2|<\dl_1$. We assume that $\dl_1<\eta$.
\par By Corollary 4.4 from \cite{P} the compact sets $K\sbs\aC$ have the Mergelyan property, i.e., there are a neighborhood $U$ of $K$ and a holomorphic mapping $h:\,U\to M$ such that $d((K,f),(K,h))<\dl_1/4$.  Let us take $\dl>0$ with the following properties: 1) any smooth Jordan neighborhood $\Om$ of $K$ that lies in the $\dl$-neighborhood of $K$ compactly belongs to $U$; 2) the diameter of $\gm\cap \Om$ is less than $\dl_1/4$; 3) the restriction of $h$ to $\ovr \Om$ that we denote also by $h$ lies in the $\dl_1/2$-neighborhood of $(K,f)$.
\par Let $\zeta$ be any point in $\gm\cap \bd\Om$. Since $(\ovr \Om,h)$ lies in the $\dl_1/2$-neighborhood of $(K,f)$ there is a point $\xi\in K$ such that $|\zeta-\xi|<\dl_1/2$ and $d(h(\zeta),f(\xi))<\dl_1/2$. Hence $|\zeta_0-\xi|<\dl_1$ because $|\zeta-\zeta_0|<\dl_1/4$. Thus $d(f(\zeta_0),f(\xi))<\eta/2$ and $d(h(\zeta),w_0)<\eta$. By Lemma \ref{L:hc}(2) we can shift $h$  so that for the shifted mapping $g$ we have $g(\zeta)=w_0$ and $(\ovr\Om,g)$ lies in the $\eps$-neighborhood of $(K,f)$.
\par If $f\in\A_{\zeta_0,w_0}(K,W,M)$ then we take $\eps>0$ so small that $\Pi^{-1}$ is defined on the $\eps$-neighborhood of $(K,\hat f)$, approximate $(K,\hat f)$ in this neighborhood and compose it with $\Pi^{-1}$.
\end{proof}
\par To continue we need the notion of {\it Rad\'o continuity.} A family of Jordan domains $\Om_t\sbs \mathbb C$, $0\le t\le1$, is called {\it Rad\'o continuous  at } $t_0\in[0,1]$  if for some $\eps>0$ a neighborhood of some point $\zeta\in\aC$ belongs to the intersection of all $\Om_t$, $t_0-\eps<t<t_0+\eps$, and the family of conformal mappings $\phi_t$ of $\mathbb D$ onto $\Om_t$ such that $\phi_t(0)=\zeta$ and $\phi_t'(0)>0$ converges uniformly on $\oD$ to $\phi_{t_0}$ as $t\to t_0$. (By a theorem of Carath\'eodory the mappings $\phi_t$  extend to $\ovr{\mathbb D}$ as its homeomorphisms onto $\ovr \Om_t$.) Such family is {\it Rad\'o continuous } if it is Rad\'o continuous  at every $t$.  A result of Rad\'o (see \cite{Ra} or \cite[Theorem II.5.2]{Go}) claims, in particular, that a family of Jordan domains $\Om_t\sbs\mathbb C$ is Rad\'o continuous if and only if for every $t_0\in[0,1]$ there are homeomorphisms $\Psi_t$ of $\bd \Om_{t_0}$ onto $\bd \Om_t$ converging uniformly to identity on $\bd \Om_{t_0}$ as $t\to t_0$.
\par The significance of Rad\'o continuity is in the following lemma.
\bL\label{L:rc} If the family of Jordan domains $\Om_t$, $t\in[0,1]$, is Rad\'o continuous, $\zeta_t\in\bd\Om_t$ is a continuous path in $\aC$ and $(\ovr\Om_t,f_t)$  is a continuous path in $\T(W,M)$ such that $f_t(\zeta_t)=w_0$, then $I_{\ovr \Om_t,\zeta_t}(f_t)\equiv\const$.
\eL
\begin{proof}
\par Suppose that $\zeta$ is a common point of all $\Om_t$ when $t$ is near $t_0$  and $\phi_t$ be the conformal mappings from the definition of Rad\'o continuity. Let $\xi_t=\phi_t^{-1}(\zeta_t)\in\aT$ and let $\al_t$ be the rotations of $\oD$ moving 1 to $\xi_t$. Since the family $\Om_t$ is Rad\'o continuous then $\xi_t$ and $\al_t$ are continuous in $t$. If $\psi_t=\phi_t\circ\al_t$ then $\psi_t(1)=\zeta_t$ and $\psi_t$ is also continuous on $\oD\times[0,1]$. Hence $(\oD,f_t\circ\psi_t)$ is a continuous path in $\A_{1,w_0}(\oD,W,M)$ and  $I_{\ovr \Om_t,\zeta_t}(f_t)\equiv\const$.
\end{proof}
\par We need the following basic lemma.
\bL\label{L:dd} Let $w_0\in W$ and $(K,f)\in\T_{\zeta_0,w_0}(W,M)$. There is $\dl>0$ such that if:
\be\item $\Om_0\sbs\sbs \Om_1$ are smooth Jordan domains;
\item pairs $(\ovr \Om_1,g_1)$ and $(\ovr \Om_0,g_0)$ lie in the $\dl$-neighborhood of $(K,f)$ in $\T(W,M)$;
\item  $\zeta_0\in\bd \Om_0$ and $\zeta_1\in\bd \Om_1$ and $g_0(\zeta_0)=g_1(\zeta_1)=w_0$;
\item there is a continuous curve $\gm:\,[0,1]\to\ovr \Om_1$  such that $\gm(t)\in \Om_1$, $0\le t<1$, $\gm(0)=\zeta_0$, $\gm(1)=\zeta_1$ and $\gm\sbs\aD(\zeta_0,\dl)$,
\ee
then $I_{\ovr \Om_0,\zeta_0}(g_0)=I_{\ovr \Om_1,\zeta_1}(g_1)$.
\eL
\begin{proof}  Let us show that if $(K,f)\in\T_{\zeta_0,w_0}(M,M)$ then for any $\eps>0$ the $\dl$ can be chosen in such a way that we can connect the  pairs $(\ovr \Om_0,g_0)$ and $(\ovr \Om_0,g_1)$ in the $\eps$-neighborhood of $(K,f)$. Let us fix $\eps>0$ and find $0<\dl_1<\eps$ such that Lemma \ref{L:hc}(1) holds with $\dl=\dl_1$. Let us denote by $\eta$ the $\dl$ in Lemma \ref{L:hc}(2)  for which this lemma holds when $\eps=\dl_1$. The mapping $f$ is uniformly continuous on $K$. So there is $\dl>0$ such that $d(f(\zeta),f(\xi))<\eta/2$ when $|\zeta-\xi|<2\dl$.  We assume that $\dl<\eta/2<\dl_1$.
\par  Since $(\ovr \Om_1,g_1)$ lies in the $\dl$-neighborhood of $(K,f)$ for any $\zeta\in\gm$ there is a point $\xi\in K$ such that $|\zeta-\xi|<\dl$ and $d(g_1(\zeta),f(\xi))<\dl$. Hence $|\xi-\zeta_0|<2\dl$ because $\gm\sbs\aD(\zeta_0,\dl)$. Thus $d(f(\zeta_0),f(\xi))<\eta/2$ and $d(g_1(\zeta),w_0)<\eta/2+\dl<\eta$.
\par Let $\Th$ be a conformal mapping of $\ovr \Om_1\sm \Om_0$ onto an annulus $A(r_0,1)=\{\zeta\in\aC:\,r_0\le|\zeta|\le1\}$ that maps $\bd \Om_1$ onto the unit circle. We define the intermediate domains $\Om_t$, $0\le t\le1$, as bounded domains with boundaries equal to $\Th^{-1}(\{|\zeta|=(1-r_0)t+r_0\})$. The domains $\Om_t$ are simply connected and the family $\Om_t$ is Rad\'o continuous because as homeomorphisms $\Psi_t$ of $\bd \Om_t$ onto $\bd \Om_{t_0}$ we can take preimages under the mapping $\Th$ of the radial correspondences between circles in $A(r_0,1)$.
\par We will reparameterize this family letting $G_t:=\Om_s$, $\gm(t)\in\bd \Om_s$, $t\in[0,1]$. Then the new family is still Rad\'o continuous. For $t\in[0,1]$ we define the pairs $(\ovr G_t,h_t)$, where $h_t$ is the restriction of $g_1$ to $\ovr G_t$. This family still lies in the $\dl$-neighborhood of $(K,f)$. Now $h_t(\gm(t))=g_1(\gm(t))$ so $d(h_t(\gm(t)),w_0)<\eta$.  By Lemma \ref{L:hc}(2) we can shift $h_t$ to get mappings $p_t$ so that $p_t(\gm(t))=w_0$ and pairs $(\ovr G_t,p_t)$ lie in the $\dl_1$-neighborhood of $(K,f)$. Note that $G_1=\Om_1$, $G_0=\Om_0$ and by the same lemma we can assume that $p_1=g_1$.
\par The pairs $(\ovr \Om_0,p_0)$ and $(\ovr \Om_0,g_0)$ are in the $\dl_1$-neighborhood of $(K,f)$ and by our choice of $\dl_1$ we can connect them by a continuous path in the intersection of the $\eps$-neighborhood of $(K,f)$ with $\A_{\zeta_0,w_0}(\ovr \Om_0,W,M)$. Consequently we can connect the pairs $(\ovr \Om_0,g_0)$ and $(\ovr \Om_1,g_1)$ in the $\eps$-neighborhood of $(K,f)$.
\par If $(K,f)\in\T_{\zeta_0,w_0}(W,M)$ then we take $\eps>0$ such that $\Pi_1^{-1}$ is defined and continuous on the $\eps$-neighborhood of $(K,\hat f)$ in $\T(M,M)$. We find $\dl>0$ for $(K,\hat f)$ such that the pairs $(\ovr \Om_0,\hat g_0)$ and $(\ovr \Om_1,\hat g_1)$ can be connected by a continuous path $(\ovr\Om_t,h_t)$, $0\le t\le 1$, in the $\eps$-neighborhood of $(K,\hat f)$ and the family of Jordan domains $\Om_t$ is Rad\'o continuous. Then the continuous path $\Pi_1^{-1}((\ovr\Om_t,h_t))$ connects  $(\ovr \Om_0,g_0)$ and $(\ovr \Om_1,g_1)$. Hence by Lemma \ref{L:rc} $I_{\ovr \Om_1,\zeta_1}(g_1)=I_{\ovr \Om_0,\zeta_0}(p_0)$.
\end{proof}
\par Now we prove that close approximations have the same homotopic type.
\bP\label{P:eht} Let $f\in\A_{\zeta_0,w_0}(K,W,M)$ and let $\gm$ be an access curve to $K$ at $\zeta_0$. There is $\dl>0$ such that if $(\ovr \Om_1,g_1)$ and $(\ovr \Om_2,g_2)$ are $\dl$-approximations of $(K,f)$ with respect to $\gm$, then $I_{\ovr \Om_2,\zeta_{\Om_2,\gm}}(g_2)=I_{\ovr \Om_1,\zeta_{\Om_1,\gm}}(g_1)$.
\eP
\begin{proof} We take as $\dl$ the $\dl$ in Proposition \ref{L:dd}. By Lemma \ref{P:phiepsa} there are a Jordan domain $\Om_0$ containing $K$ such that $\ovr \Om_0\sbs \Om_1\cap \Om_2$ and, given any point $\zeta_1\in\gm\cap\bd \Om_0$, a mapping $g_0\in\A(\ovr \Om_0,W,M)$ such that the pair $(\ovr \Om_0,g_0)$ lies in the $\dl$-neighborhood of $(K,f)$ and  $g_0(\zeta_1)=w_0$. Let $t_0=\sup\{t:\,t<s_{\Om_1,\gm},\gm(t)\in \Om_0, \}$ and let $\zeta_1=\gm(t_0)$. By Lemma \ref{L:dd} $I_{\ovr \Om_1,\zeta_{\Om_1,\gm}}(g_1)=I_{\ovr \Om_0,\zeta_1}(g_0)$ and by the same argument $I_{\ovr \Om_2,\zeta_{\Om_2,\gm}}(g_2)=I_{\ovr \Om_0,\zeta_1}(g_0)$.
\end{proof}
\par Let $\gm$ be an access curve to $K$ at $\zeta_0$. We define the mapping \[I_{K,\gm}=I_\gm:\,\H_{\zeta_0,w_0}[K,W,M]\to\H_{1,w_0}[\oD,W,M]=\eta_1(W,M,w_0)\]
as $I_{K,\gm}(f)=I_{\ovr \Om,\zeta_{\Om,\gm}}(g)$, where $(\ovr \Om,g)$ is a sufficiently close approximation of $(K,f)$. By Proposition \ref{P:eht} this mapping is well defined.
\par If $f\in\A_{w_0}(W,M)$ let $\iota(f)$ be the loop $f|_{\aT}$ in $W$. Clearly, if $[f]_{1,w_0}=[g]_{1,w_0}$ in $\eta_1(W,M,w_0)$, then $\iota(f)$ and $\iota(g)$ are homotopic in $\pi_1(W,w_0)$. Hence the mapping \[\iota_1:\,\eta_1(W,M,w_0)\to\pi_1(W,w_0)\] is also well-defined.
\par The following result shows that $I_{K,\gm}(f)$ continuously depends on $(K,f)$.
\bT\label{T:ci} For any pair $(K,f)\in\T_{\zeta_0,w_0}(W,M)$ and an access curve $\gm$ to $K$ at $\zeta_0$ there is $\dl>0$ such that $I_{K,\gm}(f)=I_{L,\gm}(g)$ for any pair $(L,g)\in\T_{\zeta_0,w_0}(W,M)$ that lies in the $\dl$-neighborhood of $(K,f)$ and has $\gm$ as an access curve to $L$ at $\zeta_0$. \eT
\begin{proof} Let $\dl_1$ be the $\dl$ from Lemma  \ref{L:dd}. Let $(\ovr \Om_1,g_1)$ be a $\dl_1$-approximation of $(K,f)$ such that the restriction of $\gm$ to $[0,s_{\ovr \Om_1,\gm}]$ lies in $\aD(\zeta_0,\dl_1)$. Let $r$ be the minimal distance between points on $\bd \Om_1$ and $K$. We take $\dl=\min\{r,\dl_1\}/2$.
\par If a pair $(L,g)$ lies in the $\dl$-neighborhood of $(K,f)$ then $L\sbs \Om_1$. We take a $\dl$-approximation $(\ovr\Om_0,g_0)$ of $(L,g)$ such that  $\Om_0\sbs\sbs \Om_1$ and $I_{\ovr\Om_0,\zeta_{\Om_0,\gm}}(g_0)=I_{L,\gm}(g)$. The pair $(\ovr\Om_0,g_0)$ lies in the $\dl_1$-neighborhood of $(K,f)$ so by Lemma \ref{L:dd}
\[I_{K,\gm}(f)=I_{\ovr\Om_1,\zeta_{\Om_1,\gm}}(g_1)=I_{\ovr\Om_0,\xi_{\Om_0,\gm}}(g_0)=I_{L,\gm}(g).\]
\end{proof}
\par The following technical lemma will be used later several times.
\bL\label{L:jd} Suppose that $K$ consists of a simple curve $\al$ connecting $\zeta_0$ and $\zeta_1$ and the closure of a smooth Jordan domain $\Om_1$ such that $\ovr \Om_1\cap \al=\{\zeta_1\}$. Let $(K,f)\in\T_{\zeta_0,w_0}(W,M)$ and let $\gm$ be an access curve to $K$ at $\zeta_0$.  Let $\Om_0\sbs \Om_1$ be another smooth Jordan domain such that $\bd\Om_1\cap\bd \Om_0=\{\zeta_1\}$ and let $L=\al\cup\ovr \Om_0$. Then there is a mapping $g\in\A_{\zeta_0,w_0}(L,W,M)$ such that $g=f$ on $\al$, $I_{\ovr \Om_1,\zeta_1}(f)=I_{\ovr \Om_0,\zeta_1}(g)$ and $I_{K,\gm}(f)=I_{L,\gm}(g)$.
\eL
\begin{proof} We take a conformal mapping $\Phi$ of $\Om_1\sm\ovr \Om_0$ onto the strip $\{0<\re\zeta<1\}$. This mapping extends smoothly to the boundary and we assume that $\Phi(\bd \Om_1)=\{\re\zeta=1\}$ and $\Phi(\bd \Om_0)=\{\re\zeta=0\}$. Since $\Phi^{-1}(\zeta)$ converges to $\zeta_1$ when $\re\zeta\to\pm\infty$ the domains $\Om_t$ bounded by curves $\Phi^{-1}(\{\re\zeta=t\})$ and $\zeta_1$ are Jordan domains. Moreover the family $\{\Om_t\}$ is Rad\'o continuous because as homeomorphisms of $\bd \Om_t$ onto $\bd \Om_{t_0}$ we can take preimages of mappings $x+it\to x+it_0$.
\par Let $\Psi_t$ be a continuous family of conformal mappings of $\Om_t$ onto $\Om_1$ such that $\Psi_t(\zeta_1)=\zeta_1$ and let $K_t=\al\cup \ovr \Om_t$. We define $f_t$ as $f$ on $\al$ and as $f\circ\Psi_t$ on $\ovr \Om_t$. Thus we obtain a continuous path in $\T_{\zeta_0,w_0}(W,M)$ and letting $g=f_0$ we get our lemma.
\end{proof}
\par Two access curves $\gm_1$ and $\gm_2$ are {\it equivalent } if for every $\eps>0$ there is $\dl>0$ such that if $0<t_1,t_2<\dl$ then the points $\gm_1(t_1)$ and $\gm_2(t_2)$ can be connected by a continuous curve $\al$ in $\mathbb D(\zeta_0,\eps)\sm K$. In the terminology of the prime ends theory (see \cite{C}) it means that curves $\gm_1$ and $\gm_2$ determine the same prime end.
\par The following result shows that $I_{K,\gm_1}=I_{K,\gm_2}$ when $\gm_1$ and $\gm_2$ are equivalent. In particular, if $K$ is the closure of a Jordan domain then by a theorem of Carath\'eodory all access curves at any point of $\bd K$ are equivalent and $I_{K,\gm}$ is determined only by $\zeta_0$ so we can write $I_{K,\gm}=I_{K,\zeta_0}$.
\bP\label{P:eht1} Let $f\in\A_{\zeta_0,w_0}(K,W,M)$ and $\zeta_0\in\bd K$. If $\gm_1$ and $\gm_2$ are equivalent access curves to $K$ at $\zeta_0$, then $I_{K,\gm_1}(f)=I_{K,\gm_2}(f)$.
\eP
\begin{proof} We take $\dl>0$ that is less than $\dl$'s in Lemma \ref{L:dd} and Propositions \ref{P:phiepsa} and \ref{P:eht}. Then we find $\dl$-approximations  $(\ovr\Om,g_1)$ of $(K,f)$ with respect to $\gm_1$ and $(\ovr\Om,g_2)$ of $(K,f)$ with respect to $\gm_2$. The domain $\Om$ has been chosen so that the restrictions of curves $\gm_1$ and $\gm_2$ to $[0,s_{\Om,\gm_1}]$ and $[0,s_{\Om,\gm_2}]$ lie in $\aD(\zeta_0,\dl)$. By Proposition \ref{P:eht} $I_{K,\gm_1}(f)=I_{\ovr \Om,\zeta_{\Om,\gm_1}}(g_1)$ and $I_{K,\gm_2}(f)=I_{\ovr \Om,\zeta_{\Om,\gm_2}}(g_2)$. We take $0<\si<\dl$ such $\aD(\zeta_0,\si)\sbs\sbs \Om$ and find $t_1,t_2>0$ such that $\gm_1(t_1)$ and $\gm_2(t_2)$ can be connected by a continuous path $\gm_3$ in $\aD(\zeta_0,\si)\sm K$.
\par We take a Jordan domain $\Om_0\sbs\sbs \Om$ such that $K\sbs \Om_0$ and $\gm_3\sbs \Om\sm\ovr \Om_0$. Let $t_3=\sup\{t:\gm_2(t)\in \Om_0\}$. By Proposition \ref{P:phiepsa} there is a $\dl$-approximation $(\ovr\Om_0,h)$ of $(K,f)$ such that $g(\gm_2(t_3))=w_0$. Let $\gm$ be a curve in $\ovr \Om\sm \Om_0$ that follows $\gm_1$ to $\gm_1(t_1)$, then $\gm_3$ until $\gm_2(t_2)$ and then $\gm_2$ to $\zeta_0$. Since $\gm_3\sbs\aD(\zeta_0,\dl)$ by Lemma \ref{L:dd} $I_{\ovr \Om,\zeta_{\Om,\gm}}(g_1)=I_{\ovr \Om_0,\gm_2(t_3)}(h)$. By the same lemma $I_{\ovr \Om,\zeta_{\Om,\gm_2}}(g)=I_{\ovr \Om_0,\gm_2(t_3)}(h)$ and we are done.
\end{proof}
\section{Holomorphic fundamental semigroup of Riemann domains}\label{S:hfs}
\par If $f\in\A_{w_0}(W,M)$ we will denote $[f]_{1,w_0}$ by $[f]$. To introduce on $\eta_1(W,M,w_0)$ a semigroup structure compatible with $\iota_1$ we need an additional construction since in the standard definition the concatenation of two loops cannot be realized as a boundary of an analytic disk.
\par Suppose that $f_1,f_2\in\A_{w_0}(W,M)$ are representatives of equivalence classes $[f_1]$ and $[f_2]$ respectively in $\eta_1(W,M,w_0)$. Let $K\sbs{\mathbb C}$ be the union of $K_1=\{|\zeta-1|\le1\}$
and $K_2=\{|\zeta+1|\le1\}$ and let $\gm(t)=-it$, $0\le t\le1$. Then $\gm$ is an access curve for $K$ to $0$. We define the mapping
\[h_{f_1,f_2}(\zeta)=\begin{cases} f_1(1-\zeta),\,\zeta\in
\bd K_1,\\ f_2(1+\zeta),\,\zeta\in \bd K_2\end{cases}\] of $\bd K$ into
$W$. The mapping $\hat h_{f_1,f_2}$ maps $K$ into $M$ so $h_{f_1,f_2}\in\A_{0,w_0}(K,W,M)$.
\par We let $[f_1]\star[f_2]=I_{K,\gm}(h_{f_1,f_2})$. If $f_1$ and $f_2$ are $h$-homotopic to $g_1$ and $g_2$ respectively in $\A_{w_0}(W,M)$, then evidently $h_{f_1,f_2}$ is $h$-homotopic to $h_{g_1,g_2}$ in $\A_{0,w_0}(K,W,M)$. Hence the class $[f_1]\star[f_2]$ is well defined.
\par The following notion of stars is similar to the notion of stars in \cite{R2}. A {\it star} is a plane compact set $K$ that consists of $n$ simple curves $\al_j:\,[0,1]\to\aC$ starting at the same point $\zeta_0$ called the {\it center} of the star and $n$ disjoint closed disks $D_j$ such that $\zeta_j=\al_j(1)\in\bd D_j$. It is required that curves $\al_j$ and $\al_i$ meet only at $\zeta_0$ when $i\ne j$ and $D_j\cap(\cup_{i=1}^n\al_i)=\{\zeta_j\}$ for all $j$. We let $K_j=\al_j\cup D_j$ and call them the {\it arms} of a star. Note that for any star there is a homeomorphism of the plane that transforms this star into a {\it straight star}, i.e., a star where all curves $\al_j$ are intervals.
\par The conformal mapping $\phi:\,\aD\to{\mathbb {CP}}^1\sm K$ extends continuously to the boundary. This happens because the natural mapping of the prime ends space of ${\mathbb {CP}}^1\sm K$ onto $\bd K$ is continuous. If $\gm$ is an access curve to $K$ at $\zeta_0$ and $\phi(1)=\zeta_0$, then the numeration of $K_j$ is chosen in such a way that as $\zeta$ travels by $\aT$ clockwise starting at $1$ the point $\phi(\zeta)$ travels first by $\bd K_1$, then $\bd K_2$ and so on.
\bP\label{P:div0} Suppose that $K$ is a star with arms $K_j=\al_j\cup D_j$, $1\le j\le n$, $\zeta_0$ is center of $K$ and $\gm$ is an access curve to $K$ at $\zeta_0$.  Let $L_2$ be the star with arms $K_j$, $2\le j\le n$, and let $L_1$ be the star with arms $K_j$, $1\le j\le n-1$. If $f\in\A_{\zeta_0,w_0}(K,W,M)$ and $f_j=f|_{K_j}$, then
\[I_{K,\gm}(f)=I_{K_1,\gm}(f_1)\star I_{L_2,\gm}(f|_{L_2})=I_{L_1,\gm}(f|_{L_1})\star I_{K_n,\gm}(f_n).\]
\eP
\begin{proof} We assume that $\zeta_0=0$. Take some small $t_0>0$ and redefine $f$ on each $\al_j$ letting it to be $w_0$ on $\al_j([0,t_0])$ while on $[t_0,1]$ we set the mappings as $f(\al_j((t-t_0)/(1-t_0)))$. We don't change $f$ on disks and preserve for new mapping the same notation. Since this operation can be achieved by a continuous family of deformations by Theorem \ref{T:ci} the $h$-homotopic classes of $(K,f)$ and all $(K_j,f_j)$ will not change.
\par In the next step for $s\in[0,1]$ we squeeze starting intervals of curves $\al_2,\dots,\al_n$ considering a continuous family of continuous curves $\al^s_2,\dots,\al^s_n$ on $[0,1]$ such that each of these curves is simple, $\al_j^0=\al_j$, $\al_j^s(t)=\al_j(t)$ for $t>t_0$ and $\al_j^1(t)=\al_i^1(t)$ for all $2\le i,j\le n$ and $t\in[0,t_1]$ for some $0<t_1<t_0$. We set $f^s_j(\al^s_j(t))=f_j(\al_j(t))$ and do not change $f_j$ on $D_j$. It is easy to see that such family can be found for straight stars and, consequently, for all stars. Let $N_2$ be the union of curves $\al^1_j$ and disks $D_j$, $2\le j\le n$, and let $g_2$ be the mapping of $N_2$ equal to $f^1_j$ on $\al^1_j\cup D_j$. Let $N=K_1\cup N_2$ and let $g\in\A_{\zeta_0,w_0}(N,W,M)$ be equal to $f_1$ on $K_1$ and to $g_2$ on $N_2$.  It follows from Theorem \ref{T:ci} that $I_{K,\gm}(f)=I_{N,\gm}(g)$ and $I_{L_2,\gm}(f)=I_{N_2,\gm}(g_2)$
\par Let $\beta_1=\al_1|_{[0,t_0]}$, $N'_1=K_1\sm\beta_1$, and $g_1=f_1|_{N'_1}$. Let $\beta_2=\al^1_2|_{[0,t_1]}$, $N'_2=N_2\sm\beta_2$, and $h_2=g_2|_{N'_2}$. We take close approximations of $(N'_j,f_j)$  by $(\ovr U_j,p_j)$, $j=1,2$, where $U_j$ are smooth Jordan domains. Proposition \ref{P:phiepsa} gives us a lot of freedom of choices for $U_j$ so we may assume
that $U_j$  meet $\beta_j$ only once at points $\zeta_j=\beta_j(s_j)$ and $\ovr U_1\cap\ovr U_2=\emptyset$. Let $\gm_j=\beta_j|_{[0,s_j]}$, $A_j=\gm_j\cup\ovr U_j$, $A=A_1\cup A_2$ and $q_j=p_j$ on $\ovr U_j$ and $w_0$ on $\gm_j$. By our choice $I_{A_j,\gm}(q_j)=I_{N_j,\gm}(g_j)$ and if approximations are close enough  $I_{A,\gm}(q)=I_{N,\gm}(g)$, where $q=q_j$ on $A_j$.
\par Using Lemma \ref{L:jd} we replace Jordan domains $U_j$ with small disjoint disks $V_j\sbs U_j$ attached to $\gm_j$ at $\zeta_j$ and the mapping $q$ on them with a mapping $r$ preserving all involved $h$-homotopy classes.  Then we take  continuous deformations  $\gm_j^s$ of $\gm_j$, $s\in[0,1]$, so that $\gm^s_j(0)=\zeta_0=0$ and $\gm^1_1(t)=t$ while $\gm^1_2(t)=-t$. We may assume that if $V_j^s=V_j+\gm_j^s(1)-\zeta_j$, $M_j^s=\gm^s_j\cup V_j^s$ and $M^s$ is the union of $M_j^s$, then $M^s$ are stars for all $s$. Let $r_j^s(\gm_j^s(t))=r_j(t)$ and $r_j^s(\zeta)=r_j(\zeta-\gm_j^s(1)+\zeta_j)$ on $V^s_j$. We let $r^s$ to be equal to $q_j^s$ on $M^s_j$.
\par Applying a continuous family of rotation and dilations we can make the disks $V^1_j$ perpendicular to the real axis and of radius 1. The mapping $r^1$ will follow these changes as compositions with rotations and dilations.
\par The obtained compact set $L$ that consist of intervals $[0,1]$ and $[-1,0]$ and disks $\oD(2,1)$ and $\oD(-2,1)$ has two prime ends at $0$: one of them is equivalent to $[-i,0]$ and another to $[0,i]$. The chosen numeration tells us that our access curve is $\gm=[-i,0]$. Now we shrink intervals $[0,1]$ and $[-1,0]$ to 0 simultaneously translating disks $V^1_j$. We can do this because $r^1_j\equiv w_0$ on these intervals. As the result we obtain figure form the definition of the $\star$ operation and we see that $I_{K,\gm}(f)=I_{K_1,\gm}(f_1)\star I_{L_2,\gm}(f|_{L_1})$.
\par The second equality in the proposition is proved similarly.
\end{proof}
\par This proposition leads to the following theorem.
\bT\label{T:sg} The operation $\star$ induces on $\eta_1(W,M,w_0)$ the structure of a semigroup with unity.
\eT
\begin{proof}  The unity is the class of the constant mapping equal to $w_0$ on $\aT$. If, say, in the definition of the $\star$ operation $f_1\equiv w_0$ then continuously shrinking $K_1$ to the origin leaving the functions equal to $w_0$ we will get a continuous path in $\T_{0,w_0}(W,M)$ which ends at $(K_2,f_2(1+\zeta))$. By Theorem \ref{T:ci} $I_\gm(h_{f_1,f_2})=[f_2]$.
\par To prove that the operation $\star$ is associative we consider a compact set $L$ consisting of three intervals $I_1=[0,1]$, $I_2=[0,i]$, $I_3=[-1,0]$ and three closed disks $D_1=\{|\zeta-2|\le1\}$, $D_2=\{|\zeta-2i|\le1\}$ and $D_3=\{|\zeta+2|\le1\}$. The access curve $\gm=[-i,0]$. Given $f_1,f_2,f_3\in\A_{1,w_0}(\oD,W,M)$ we define the mapping $f$ on $L$ to be equal to $w_0$ on intervals $I_1, I_2, I_3$ and $f_1(2-\zeta)$ on $D_1$, $f_2(2+i\zeta)$ on $D_2$ and $f_3(2+\zeta)$ on $D_3$.
\par By Proposition \ref{P:div0}
\[I_{L,\gm}(f)=I_{L_1,\gm}(f_1)\star (I_{L_2,\gm}(f_2)\star I_{L_3,\gm}(f_3))=(I_{L_1,\gm}(f_1)\star I_{L_2,\gm}(f_2))\star I_{L_3,\gm}(f_3).\]
\end{proof}
\par Using induction and Proposition \ref{P:div0} and the previous theorem we get
\bT\label{T:div} In assumptions of Proposition \ref{P:div0}
\[I_{K,\gm}(f)=I_{K_1,\gm}(f_1)\star\cdots\star I_{K_n,\gm}(f_n).\]
\eT
\par We finish this section with two examples of the semigroup $\eta_1$ when $W=A_{s,r}=\{ s< |z| <r \}$, where $0<s<1<r$, and $M={\mathbb {CP}}^1$ or $M=\aD(0,x)$, where $r\le x\le\infty$. We fix $\Pi(z)=z$ and $w_0=1$. The examples below show that the mapping $\iota_1:\,\eta_1(W,M,w_0)\to\pi_1(W,w_0)$ need not to be injective or surjective.
\par If $f\in\A_{w_0}(A_{s,r},{\mathbb {CP}}^1)$ then we can write it as $\hat f=gB_1B_2^{-1}$, $B_1$ and $B_2$ are Blaschke products and $B_1$ contains all zeros of $\hat f$ while $B_2$ contains all poles. The function $g$ has no zeros and poles, $g(1)=1$ and $s<|g|<r$ on $\aT$. Hence $g$ maps $\oD$ into $A_{s,r}$ and is $h$-homotopic to the constant mapping equal to 1.
\par We change $B_1$ by dragging its zeros to 0 by continuous curves and then change $B_2$ by dragging its zeros to some $a\ne0$ at $\aD$. Thus $f$ is $h$-homotopic to
\[h(\zeta)=(-a)^n\zeta^m\left(\frac{1-\bar a\zeta}{\zeta-a}\right)^n\]
Thus we obtained a mapping $[f]\to(m,n)$ and, clearly, it is a homomorphism and it is injective because continuous deformations of analytic disks do not change the numbers of zeros and poles. Hence the semigroup $\eta_1(A_{s,r},{\mathbb {CP}}^1, w_0)$ is isomorphic to $\mathbb N_0\oplus\mathbb N_0$, where $\mathbb N_0$ is the semigroup by addition of non-negative integers.
\par A similar argument shows that the semigroup $\eta_1(A_{s,r},\aD(0,x), w_0)$ is isomorphic to $\mathbb N_0$. The  isomorphism is given by the mapping $[f]\to m$, where $m$ is the number of zeros of $\hat f$ counted with multiplicity.
\section{Properties of holomorphic fundamental semigroups}\label{S:phf}
\par Let $(W_1,\Pi_1)$ and $(W_2,\Pi_2)$ be two Riemann domains over two complex
manifolds $M_1$ and $M_2$ respectively. Suppose $w_1\in W_1$, $w_2\in W_2$ and there are holomorphic mappings $\phi: W_1 \to W_2$ such that $\phi(w_1)=w_2$ and $\psi: M_1 \to M_2$ which satisfy $\psi\circ \Pi_1=\Pi_2\circ \phi$. Then for any $f \in \A(K,W_1,M_1)$ we have $\Pi_2\circ\phi \circ f=\psi\circ \Pi_1\circ f=\psi \circ \hat f$. So
$\widehat{\phi\circ f}=\psi\circ \hat f$ and we get a continuous mapping from
$\T(W_1,M_1)$ to $\T(W_2,M_2)$ which maps a pair $(K,f)$ to $(K,\phi\circ f)$. Hence, firstly, the mapping from $\A_{w_1}(W_1,M_1)$ to $\A_{w_2}(W_2,M_2)$ induces a well defined mapping $\phi_*$ from $\eta_1(W_1,M_1,w_1)$ to $\eta_1(W_2,M_2,w_2)$ given by $\phi_*([f])=[\phi \circ f]$. Secondly, if $\gm$ is an access curve to $K$, then $\phi_*(I_{K,\gm}(f))=I_{K,\gm}(\phi\circ f)$. In particular, if $(K,h_{f_1,f_2})$ is the pair in the definition of the $\star$ operation, then
\begin{equation}\begin{aligned}&\phi_*([f_1]\star[f_2])=\phi_*(I_{K,\gm}(h_{f_1,f_2}))\notag\\&=I_{K,\gm}(\phi\circ h_{f_1,f_2})=I_{K,\gm}(h_{\phi\circ f_1,\phi\circ f_2})=[\phi\circ f_1]\star[\phi\circ f_2]=\phi_*[f_1]\star\phi_*[f_2].
\end{aligned}\end{equation}
\par This leads us to the following proposition.
\bP\label{P:mo} The induced mapping $\phi_*:\eta_1(W_1,M_1,w_1) \to \eta_1(W_2,M_2,w_2)$ is a homomorphism.
\eP
\par Clearly $(W_1 \times W_2, (\Pi_1,\Pi_2))$
is a Riemann domain over $M_1\times M_2$. As in the classical homotopy theory Proposition \ref{P:mo} leads to the following corollary.
\bC If $(W_1,\Pi_1)$ and $(W_2,\Pi_2)$ are two Riemann domains over two complex
manifolds $M_1$ and $M_2$ respectively, then
\[  \eta_1(W_1\times W_2,M_1 \times M_2, (w_1,w_2)) \cong \eta_1(W_1,M_1,w_1)
\times \eta_1(W_2,M_2,w_2).  \]
\eC
\par Another corollary describes the powers in $\eta_1(W,M,w_0)$.
\bC Let $f\in\A_{w_0}(W,M)$ and let $[f]^{\star k}$ be the product of $k$ classes $[f]$. Then $[f]^{\star k}=[f(\zeta^k)]$.
\eC
\begin{proof} We may assume that $\hat f$ is defined on $\aD(0,r)$, $r>1$, and $f$ maps $A_{r^{-1},r}=\{\zeta\in\aC:r^{-1}<|\zeta|<r\}$ into $W$. Set $W_1=A_{r^{-1},r}$, $M_1=\aD(0,r)$ and $w_1=1$. Let $\phi=f$. By Proposition \ref{P:mo} and an example at the end of the previous section we have
\[[f(\zeta^k)]=\phi_*([\zeta^k])=\phi_*([\zeta]^{\star k})=\phi_*([\zeta])^{\star k}=[f]^{\star k}.\]
\end{proof}
\par Let $\al(t)$, $t \in [0,1]$, be a continuous curve in $W$ with $\al(0)=w_0$ and $\al(1)=w_1$. Let $L$ be a compact set on the plane consisting of the interval $I=[0,1]$ and the disk $D=\{|\zeta-2|\le1\}$. Given a mapping $f\in\A_{1,w_1}(\oD,W,M)$ we define a mapping $\wtl f$ on $L$ to be equal to $\al$ on $I$ and to $f(2-\zeta)$ on $\bd D$. Clearly, $\wtl f\in\A_{0,w_0}(L,W,M)$.
\par We take the access curve $\gm(t)=-it$, $0\le t\le 1$, to $L$ at the origin. Clearly, if $[f]_{1,w_0}=[g]_{1,w_0}$, then $I_{L,\gm}(\wtl f)=I_{L,\gm}(\wtl g)$. Hence we have a well-defined mapping $F_\al(f)=F_\al([f])=I_{L,\gm}(\wtl f)$ from $\eta_1(W,M,w_1)$ into $\eta_1(W,M,w_0)$.
\par By Theorem \ref{T:ci} any curve connecting $w_0$ to $w_1$ which is homotopic to $\al$ will give us the same mapping $F_{\al}$. Thus $F_\al$ depends only on the homotopy class $\{\al\}$ of $\al$ in $\pi_1(W,w_0,w_1)$.
\par We let $\al^{-1}$ to be the curve $(\al^{-1})(t)=\al(1-t)$ for $0 \le t \le 1$ and, if a curve $\beta$ connects $w_1$ and $w_2\in W$, denote by $\al\beta$ the curve on $[0,1]$ defined as $\al(2t)$ when $0\le t\le 1/2$ and as $\beta(2t-1)$ when $1/2\le t\le1$.
\bT\label{T:bp} Let $w_0, w_1,w_2$ are points of $W$, continuous curves $\al$ and $\beta$ connect $w_0$ with $w_1$ and $w_1$ with $w_2$ respectively. Then:
\be\item $F_{\al\beta}=F_{\al}\circ F_{\beta}$;
\item $F_\al$ is an isomorphism of $\eta_1(W,M,w_1)$ onto $\eta_1(W,M,w_0)$.
\ee
\eT
\begin{proof} (1) Let $K=[0,2]\cup\oD(3,1)$ and let $f\in\A_{w_2}(W,M)$. We define the mapping $g$ on $K$ as $\al$ on $[0,1]$, $\beta(t-1)$ on $[1,2]$ and $f(3-\zeta)$ on $\oD(3,1)$. All access curves to $K$ at 0 are equivalent, we take any such $\gm$ and $I_{K,\gm}(g)=F_{\al\beta}(f)$.
\par We take a Jordan domain $\Om$ containing $K_1=[1,2]\cup\oD(3,1)$ in its closure such that $1\in\bd \Om$ and a close approximation $h_1$ of $g|_{K_1}$ on $\ovr \Om$ so that $I_{\ovr \Om,1}(h_1)=I_{K_1,1}(g|_{K_1})=F_\beta(f)$. If the mapping $h$ is defined on $K_2=[0,1]\cup\ovr \Om$ as $h_1$ on $\ovr \Om$ and $\al$ on $[0,1]$, then also $I_{K,\gm}(f)=I_{K_2,0}(h)$.
\par Then using Lemma \ref{L:jd} we deform $\Om$ to a small disk attached to 1 and then to the disk $\aD(2,1)$ while $h$ changes to $p_1$. If the mapping $p$ is defined on $K_3=[0,1]\cup\oD(2,1)$ as $p_1$ on $\oD(2,1)$ and as $\al$ on $[0,1]$, then also $I_{K,\gm}(f)=I_{K_3,\gm}(p)$. Since $I_{\oD(2,1),1}(p)=I_{\ovr \Om,1}(h_1)$ we see that $I_{K_3,\gm}(p)=F_\al(p_1)=F_\al(F_\beta((f))$.
\par (2) Let us show that $F_\al$ is a homomorphism, i.e., if  $[f]=[f_1]\star[f_2]$ then $F_\al(f)=F_\al(f_1)\star F_\al(f_2)$. Let $K$ be the compact set from the definition of the $\star$ operation, $K'=[-i,0]\cup K$ and the mapping $h$ is equal to $h_{f_1,f_2}$ on $K$ and as  $\al(-it+1)$ on $[-i,0]$.
\par We take the interval $[-2i,-i]$ as an access curve $\gm$ to  $K'$ at $-i$ and consider compact sets $N_1=[-i,0]\cup K_1$ and $N_2=[-i,0]\cup K_2$. Then we rotate $N_1$ by a small positive angle around $-i$ and rotate $N_2$ by a small negative angle around $-i$. The mapping $h$ follows these rotations. The obtained set is a star $(P,q)$ with two arms $(P_1,q_1)$ and $(P_2,q_2)$ such that $I_{P_1,\gm}(q_1)=F_\al(f_1)$ and $I_{P_2,\gm}(q_2)=F_\al(f_2)$. By Theorem \ref{T:ci} $I_{P,\gm}(q)=I_{K',\gm}(h)=F_\al(f)$ and by Theorem \ref{T:div} \[I_{P,\gm}(q)=I_{P_1,\gm}(q_1)\star I_{P_2,\gm}(q_2)=F_\al(f_1)\star F_\al(f_2).\] Thus $F_\al$ is a homomorphism. Since by (1) $F_{\al^{-1}}\circ F_{\al}$ is an identity mapping, $F_\al$ is an isomorphism.
\end{proof}
\par If $g\in\A_{w_0}(W,M)$ then we let the mapping $F_g=F_\al$, where $\al(t)=g(e^{2\pi it})$, and if $\al$ is a loop in $W$ starting at $w_0$ we denote by $\{\al\}$ the equivalence class of $\al$ in $\pi_1(W,w_0)$.
\bT\label{T:ap} The mapping $\Phi:\,\{\al\}\to F_\al$ establishes a homomorphism of $\pi_1(W,w_0)$ into $\Aut(\eta_1(W,M,w_0))$. Moreover, if $g\in\A_{w_0}(\oD, W,M)$, then $F_{g}([f])\star [g]=[g]\star[f]$ and $F_g([g])=[g]$.
\eT
\begin{proof} The first part of the theorem is a direct consequence of Theorem \ref{T:bp}. Let us show that $F_{g}([f])\star [g]=[g]\star[f]$.
\par Consider a compact set $K$ consisting of the disks $D^1=\{|\zeta|\le1\}$ and $D^2=\{|\zeta-3|\le1\}$ and the interval $I=[1,2]$. We define a mapping $h(\zeta)$ on $D^1$ as $g(\zeta)$  and on $D^2$ as $f(3-\zeta)$. Let $h(t)=g(e^{2\pi it})$ on $I$. Let $\gm=[-i+1,1]$ be an access curve to $K$ at $1$. Then $I_{K,\gm}(h)=F_g([f])\star[g]$.
\par Consider the continuous family of compact sets $K_s$, $0\le s\le 1/2$, consisting of the disk $D^1$, an interval $I_s=[e^{2\pi si},(2-s)e^{2\pi si}]$ and the closed unit disk $D^2_s$ attached normally to $(2-s)e^{2\pi si}$. The mapping $h_s$ on $K_s$ is defined as $h$ on $D^1$ and as $h(s+|\zeta|)$ when $\zeta\in I_s$. The mapping $h_s$ on $D^2_s$ is defined as a composition of $h$ on $D^2$ and a conformal mapping that maps $D^2_s$ onto $D^2$ moving $(2-s)e^{2\pi si}$. Simply speaking we rotate $I\cup D^2$ around $D^1$ leaving one end of $I_s$ attached normally to $D^1$. Clearly, the pairs $(K_s,h_s)$ form a continuous path and $I_{K_s,\gm}(h_s)=F_g([f])\star[g]$.
\par When $s=1/2$ the set $K_{1/2}$ consists of $D^1$, $I_{1/2}=[-1,-3/2]$ and the disk $D^2_{1/2}$. Since all access curves to $K_{1/2}$ at $1$ are equivalent we replace $\gm$ with $\gm'=[i+1,1]$. Still $I_{K_{1/2},\gm'}(h_{1/2})=F_g([f])\star[g]$.
\par Then we continue the process described above for $1/2\le s\le1$. Finally, $K_1$ will consists of $D^1$ and $D^2_1=\{|\zeta-2|\le1\}$. The mapping $h_1$ is equal to $g$ on $D^1$ and to $f(2-\zeta)$ on $D^2_1$. Now it is clear that $I_{K_1,\gm'}(h_1)=[g]\star [f]$.
\par To show that $F_g([g])=[g]$ we start with the compact set $K_1$ consisting of the interval $I=[0,1]$ and the unit disk $D_1=\{|\zeta-2|\le1\}$. The mapping $f_1$ on $K_1$ is defined as $g(e^{2\pi it})$ on $I$ and as $g(2-\zeta)$ on $D_1$. If the access curve $\gm=[-i,0]$, then $I_{K_1,\gm}(f_1)=F_g([g])$.
\par  For $0\le s\le1$ we define compact sets $K_s$ consisting of the intervals $I_s=[0,s]$ and the disks $D_s=\{|\zeta-(1+s)|\le1\}$. The mapping $f_s$ is defined as $g(e^{2\pi it})$ on $I_s$ and as $g(e^{2\pi is}(1+s-\zeta))$ on $D_s$. The pairs $(K_s,f_s)$ form a continuous path and $I_{K_s,\gm}(f_s)=F_g([g])$. Since $K_0$ consists of the disk $\{|\zeta-1|\le1\}$ and the mapping $f_0(\zeta)=g(1-\zeta)$ we see that $I_{K_1,\gm}(f_1)=[g]$.
\end{proof}
\par We remind that a semigroup $S$ is {\it cancellative} if $ab=ac$ or $ba=ca$ imply $b=c$ for any $a,b,c\in S$ and it is {\it right(left)-reversible} if $Sa\cap Sb\ne\emptyset$ ($aS\cap bS$) for any $a,b\in S$ and {\it reversible} if it is both right- and left- reversible.
\bC\label{C:cor} Let $f,g\in\A_{w_0}(W,M)$. Then: \be\item if $\iota_1([f])=\iota_1([g])$ then $[f]\star[g]=[g]\star[f]$.
\item The semigroup $\eta_1$ is reversible.
\item The semigroup $\eta_1$ is embeddable into a group if and only if it is cancellative.
\item The image of $\eta_1(W,M,w_0)$ in $\pi_1(W,w_0)$ under the mapping $\iota_1$ is invariant with respect to the inner automorphisms in $\pi_1(W,w_0)$.
\ee
\eC
\begin{proof} Since $F_f=F_g$, $[f]\star[g]=[g]\star[f]$ and we get (1) by the second part of Theorem \ref{T:ap}. (2) follows because $F_{g}([f])\star [g]=[g]\star[f]$ and $[f]\star[g]=[g]\star F_{g^{-1}}([f])$. (3) follows from Ore's theorem (\cite[1.10]{CP}) which says that any right-reversible  cancellative semigroup can be embedded into a group.
\par To show that the image is invariant with respect to the inner automorphisms we take $f\in\A_{w_0}(W,M)$ and the representative $\al$ of $\{\al\}\in\pi_1(W,w_0)$. Since $\{\al\}\{f\}\{\al\}^{-1}=\iota_1(F_\al([f]))$ the invariance follows.
\end{proof}
\par We say that elements $[f_0]$ and $[f_1]$ in $\eta_1(W,M,w_0)$ are {\it $\pi_1$-conjugate} if $F_{\al}([f_0])=[f_1]$ for some $\al\in\pi_1(W,w_0)$. The sets of all elements of $\eta_1(W,M,w_0)$ that are $\pi_1$-congugate of $[f]$ is said to be the {\it $\pi_1$-conjugacy class} of $[f]$. We denote the set of all $\pi_1$-conjugacy classes by $\C(W,M,w_0)$.
\bP\label{P:fh} Mappings $f_0,f_1\in\A_{w_0}(W,M)$ belong to the same connected component of $\A(W,M)$ if and only if they are $\pi_1$-conjugate.
\eP
\begin{proof} If mappings $f_0,f_1\in\A_{w_0}(W,M)$ belong to the same connected component of $\A(W,M)$ then there is a continuous curve $f_t$, $0\le t\le 1$, in $\A(W,M)$  that connects $f_0$ and $f_1$. For $0\le t \le 1$ define $\al(t)=f_t(1)$, $K_t=[0,t]\cup\oD(1+t,1)$ and the mappings $g_t$ to be equal to $\al$ on $[0,t]$ and to $f_t(1+t-\zeta)$ on $\oD(1+t,1)$. The access curve $\gm$ is $[-1,0]$. Note that the curve $\al$ is closed, $I_{K_1,\gm}([g_1])=F_\al([f_1])$ and $I_{K_0,\gm}([g_0])=[f_0]$. By Theorem \ref{T:ci} $F_\al([f_1])=[f_0]$.
\par If $F_{\al}([f_1])=[f_0]$ for some $\al\in\pi_1(W,w_0)$ then we let $K=[0,1]\cup\oD(2,1)$ and let $g$ to be equal to $\al$ on $[0,1]$ and to $f_1(2-\zeta)$ on $\oD(2,1)$. The access curve $\gm$ is $[-i,0]$. Then $I_{K,\gm}(g)=[f_0]$ and $I_{\oD(2,1),1}(g)=[f_1]$.
\par Let $(\ovr\Om,h)$ be a close approximation of $(K,g)$ with the following properties: $\ovr \Om$ is symmetric with respect to the real axis and intersects the real axis only at 0 and $x_0>3$, $h(0)=h(1)=w_0$, $I_{\Om,0}(h)=[f_0]$ and $I_{\oD(2,1),1}(h)=[f_1]$. We take a conformal mapping $\Phi$ of $\Om\sm \oD(2,1)$ onto an annulus $A=\{s\le|\zeta|\le r\}$ such that $\Phi(\bar z)=\bar\Phi(z)$ and  $\Phi$ moves $[0,1]$ to $[s,r]$. Define domains $\Om_t$ as domains bounded by $\Phi^{-1}(\aT(0,t))$, $s\le t\le r$. The domains $\Om_t$ are Jordan and if $e_t$ are conformal mapping of $\oD$ onto $\Om_t$ such that $e_t(1)=\Phi^{-1}(t)$ and $e'_t(0)>0$, then the mappings $h\circ e_t$ form a continuous path in $\A(W,M)$ while $[h\circ e_r]=[f_1]$ and $[h\circ e_s]=[f_0]$.
\end{proof}
\par It follows from this proposition that the mapping $\Psi$ assigning to each class in $\C(W,M,w_0)$ the connected component of $\A(W,M)$ containing representatives of elements in this class is well defined.
\bT The mapping $\Psi$ is a bijection.
\eT
\begin{proof} By Proposition \ref{P:fh} $\Psi$ is an injection. If $U$ is a connected component of $\A(W,M)$ then there are a point $w_1\in W$ and a mapping $f\in\A_{w_1}(W,M)$ such that $f\in U$. Let $\al$ be a path in $W$ that connects $w_0$ to $w_1$. By Theorem \ref{T:bp} $F_{\al}([f])\in\eta_1(W,M,w_0)$ and the same argument as in the proof of ``only if'' part in Proposition \ref{P:fh} shows that the representatives of $F_{\al}([f])$ are in $U$. So $\Psi$ is surjective.
\end{proof}
\section{The group $\rho_1(W,M,w_0)$}\label{S:rho}
\par By the analogy with the complex case we introduce the space $\T^{\aR}(W,M)$ of pairs $(K,f)$, where $K$ is a connected compact set on the plane with connected complement and $f$ is a continuous mapping of $\bd K$ into $W$ such that $\hat f=\Pi\circ f$ extends to a continuous mapping of $K$ into $M$. If $(K,f)$ and $(L,g)$ are in $\T^{\aR}(W,M)$ we define the distance between  $(K,f)$ and $(L,g)$ similar to the definition of the distance $d$ on $\T(W,M)$. That makes the imbedding of $\T(W,M)$ into $\T^{\aR}(W,M)$ an isometry.
\par  For a such compact set $K$, a point $\zeta_0\in\bd K$ and a point $w_0\in W$ let us denote by $\R_{\zeta_0,w_0}(K,W,M)$ the subset of all pairs $(K,f)\in\T(W,M)$ such that $f(\zeta_0)=w_0$. If $K=\oD$ then $\R_{w_0}(W,M)=\R_{1,w_0}(\oD,W,M)$. We say that $f_0,f_1\in \R_{\zeta_0,w_0}(K,W,M)$ are equivalent if they belong to the same connected component of $\R_{\zeta_0,w_0}(K,W,M)$ and denote the set of all equivalence classes by $\H^{\aR}_{\zeta_0,w_0}(K,W,M)$ and let $[f]_\rho$ be the equivalence class containing $f$.
\par If $\gm$ is an access curve to $K$ at $\zeta_0$ then, similar to the complex case, we can introduce the mapping
\[I^{\aR}_{K,\gm}:\,\H^{\aR}_{\zeta_0,w_0}(K,W,M)\to\H^{\aR}_{1,w_0}(\oD,W,M)=\rho_1(W,M,w_0).\]
Similar to the $\star$ operation introduced earlier we can define the $\star$ operation on $\rho_1(W,M,w_0)$ and a similar but simpler reasoning shows that $\rho_1(W,M,w_0)$ with the $\star$ operation is a semigroup with unity. All properties of the operator $I_{K,\gm}$ and the semigroup $\eta_1$, proved in the previous sections, stay true for their analogs $I^{\aR}_{K,\gm}$ and $\rho_1$ but $\rho_1$ is a group.
\bT\label{T:rho}  The operation $\star$ induces on $\rho_1(W,M,w_0)$ the structure of a group: if $[f]_\rho\in\rho_1(W,M,w_0)$ then $[f(\bar\zeta)]_\rho=[f]^{-1}_\rho$.\eT
\begin{proof}
\par Let $f_1\in\R_{w_0}(W,M)$ and let $f_2(\zeta)=f_1(\bar\zeta)$. Let  $K_1=\{|\zeta-1|\le1\}$ and $K_2=\{|\zeta+1|\le1\}$ be the sets from the definition of the $\star$ operation. Let $K_1^t=K_1-t$ and $K^t_2=K_2+t$ and let $L_t=K^t_1\cup K^t_2$ when $0\le t\le 1$ and $L_t=K_1^t\cap K_2^t$ when $1\le t\le2$. Since  $\hat f_2(1+\zeta)=\hat f_1(1+\bar\zeta)=\hat f_1(1-(-\bar\zeta))$ the mappings $\hat f_1(1-\zeta)$ and $\hat f_2(1+\zeta)$ are symmetric with respect to the imaginary axis. Hence the mappings $g_t$ equal $\hat f_1(1-\zeta+t)$ on $\bd L_t\cap K_1^t$ and $\hat f_2(1+\bar\zeta-t)$ on $\bd L_t\cap K_2^t$ are continuous and $\hat g_t$ extends to the continuous mapping of $L_t$.
\par To preserve the base points we shift $L_t$ and $f_t$ upward by $i\al(t)$,  where $\al(t)=\cos^{-1}(1-t)$, and let $L'_t=(L_t+i\al(t))\cup[0,i\al(t)]$. We set $g'_t$ as $g_t$ shifted upward and also let $g'_t(iy)=f_1(e^{i\al(y)})$ for $0\le y\le\al(t)$. Note $g_t(0)=f_1(1)=w_0$ and $g_t(i\al(t))=f_1(e^{i\al(t)})$.
\par The path $(L'_t,g'_t)$, $0\le t\le2$, is continuous in $\T^{\aR}(W,M)$ and by the real analog of Theorem \ref{T:ci} we see that $[f_1]_\rho\star[f_2]_\rho=[g'_2]_\rho$. But $L'_2$ is the interval $[0,i\pi]$ so $[g'_2]_\rho=e$.
\end{proof}
\par There are natural homomorphisms $\dl_1:\,\eta_1(W,M,w_0)\to\rho_1(W,M,w_0)$ and $\dl_2:\,\rho_1(W,M,w_0)\to\pi_1(W,w_0)$ such that $\dl_2\circ\dl_1=\iota_1$.
\bT\label{T:dls} Let $W\sbs M$. Then in notation above:
\be\item if $\pi_1(M,w_0)=0$ then $\dl_2$ is onto;
\item if $\pi_2(M,w_0)=0$ then $\dl_2$ one-to-one.
\ee
\eT
\begin{proof} (1) is evident. To show (2) we take an element $[f]_\rho\in\ker\dl_2$ and let $\al=f|_\aT$. Then $\{\al\}=\dl_2([f]_\rho)$. There is a continuous mapping $g:\,\oD\to W$ such that $g|_{\aT}=\al$. It means that $[g]_\rho=e$. If we realize $\hat f$ as a mapping of the upper hemisphere of the unit ball in $\aR^3$ and $g$ as the mapping of the lower one then we obtain the mapping $h$ of the sphere $S^2$ into $M$ equal to $\al$ on the equator. We may assume that $h(1,0,0)=w_0$. Since $\pi_2(M,w_0)=0$ the mapping $h$ can be continuously extended to the ball as a mapping into $M$. Thus $[f]_\rho=[g]_\rho=e$.
\end{proof}
\par As simple consequences of Corollary \ref{C:cor}(1) and Theorem \ref{T:dls} we obtain
\bC Let $W\sbs M$. Then \be\item The kernel of $\dl_1$ is a commutative semigroup.
\item If $\pi_1(M,w_0)=\pi_2(M,w_0)=0$ then $\rho_1(W,M,w_0)=\pi_1(W,w_0)$.
\ee
\eC
\par If $W\sbs M$ and $\Pi$ is an inclusion map then $\rho_1(W,M,w_0)$ is, of course, the relative homotopy group $\pi_2(M,W,w_0)$. So for examples in Section \ref{S:hfs} we get that
$\rho_1(A_{s,r},\mathbb{CP}^1,w_0)=\aZ\oplus\aZ$ while $\rho_1(A_{s,r},\aC,w_0)=\aZ$.
\section{Complements to analytic varieties}\label{S:ctav}
\par Let $A$ be an analytic set in a connected complex manifold $M$ and $W=M\sm A$. We assume that $A$ is is the union of irreducible components $A_j$ of pure codimension $1$. (Analytic sets of codimension $2$ and higher do not influence groups $\eta_1$ and $\rho_1$.)
\par  If $f\in\R_{w_0}(W,M)$ and the set $A(f)=\{\zeta\in\aD:\, f(\zeta)\in A\}$ is finite, then we define the index $\ind(f,A_j)$ as the intersection index of $f(\oD)$ and $A_j$. If $\zeta\in A_j(f)$ and $\phi$ is a defining function of $A_j$ on a neighborhood of $f(\zeta)$, then we define a local index $\ind_\zeta(f,A_j)$ of $f$ at $\zeta$ as the index of $\phi\circ f$ at $\zeta$. Hence
\[\ind(f,A_j)=\sum_{f(\zeta_k)\in A_j}\ind_{\zeta_k}(f,A_j).\]
A general $f\in\R_{w_0}(W,M)$ can be approximated by such mappings and close approximations have the same indexes so $\ind(f,A_i)$ is defined for all $f\in\R_{w_0}(W,M)$. The index is a homotopic invariant so if $f_0,f_1\in\R_{w_0}(W,M)$ and $[f_0]_\rho=[f_1]_\rho$, then $\ind(f_0,A_j)=\ind(f_1,A_j)$. Thus the mapping $\ind$ is well defined on $\rho_1$. Also $\ind(f,A_j)=0$ for all $j$ if $[f]_\rho=e$.
\par It follows directly from the definition of the $\star$ operation that $\ind([f]_\rho\star[g]_\rho,A_j)=\ind([f]_\rho,A_j)+\ind([g]_\rho,A_j)$. In particular, the group $\rho_1$ has no idempotents.
\par Suppose that $f\in\R_{w_0}(W,M)$ and the set $A(f)$ is finite. Let $K$ be a star in $\oD$ with its center at $1$ such that its arms $K_j$ consist of simple curves $\al_j\in\oD\sm A(f)$ that meet $\bd\oD$ only at 1 and closed disjoint disks $D_j\sbs\aD$, $1\le j\le k$, such that the set $\bd D_j\cap A(f)$ is empty, each $D_j$ contains exactly one point of $A(f)$ and $A(f)$ is covered by disks $D_j$. Let $\wtl f$ be the restriction of $f$ to $K$. We will call $(K,\wtl f)$ the {\it factorization } of $f$.
\bT\label{T:ft} Suppose that $f\in\R_{w_0}(W,M)$ and $f(\aD)$ meets $A$ only at finitely many points $\zeta_1,\dots,\zeta_k$. Let $K=\cup_{j=1}^kK_j$ be a factorization of $f$. Then
\[[f]_\rho=I^{\aR}_{K,1}(\wtl f)=\prod_{j=1}^kI^{\aR}_{K_j,1}(\wtl f).\]
\eT
\begin{proof} Since in this case $\rho_1(W,M,w_0)=\pi_2(M,W,w_0)$ and $K$ is a homotopic retract of $\oD$, the first equality follows. The second equality follows from the real analog of Theorem \ref{T:div}.
\end{proof}
\par Now we can prove the general analog of the result of L. Rudolph in \cite{R1}. Let $S=\dl_1(\eta_1(W,M,w_0))$ and $S^{-1}$ is the semigroup consisting of all $a\in\rho_1(W,M,w_0)$ such that $a=b^{-1}$ and $b\in S$.
\bT\label{T:rr} The semigroup $S$ has the following properties:
\be\item $S$ is reversible;
\item $S\cap S^{-1}=\{e\}$;
\item any element $a\in\rho_1(W,M,w_0)$ is expressible in the form $bc^{-1}$, $b,c\in S$;
\item any element $a\in\rho_1(W,M,w_0)$ is expressible in the form $d^{-1}f$, $d,f\in S$.
\ee
\eT
\begin{proof} (1) holds because $\eta_1(W,M,w_0)$ is reversible. To show (2) we suppose that $a=b^{-1}$, $a=\dl_1([f_0])\ne e$, and $b=\dl_1([f_1])$, $f_0,f_1\in\A_{w_0}(W,M)$. If the set $A(f_0)$ is empty then $[f_0]=e$ and $[f_1]=e$. So we assume that $\ind(f_0,A_j)>0$ for some $j$. Let $f_2(\zeta)=f_1(\bar\zeta)$. By Theorem \ref{T:rho} $[f_2]_\rho=[f_1]_\rho^{-1}=[f_0]_\rho$. But $\ind(f_0,A_j)>0$ while $\ind(f_2,A_j)<0$ and we came to a contradiction.
\par To show (3) we take $f\in\R_{w_0}(W,M)$ that is smooth and transverse to $A$ and $[f]_\rho=a$. In this case the set $A(f)$ is finite and consists of points $\zeta_1,\dots,\zeta_k$ in $\aD$ such that $\ind_{\zeta_j}(f,A)=\pm1$. We assume that points $\zeta_j$ are enumerated in such a way the this local index is 1 when $1\le j\le n$ and $-1$ when $n<j\le k$ and change $f$ slightly near these points so it become holomorphic when $1\le j\le n$ and antiholomorphic when $n<j\le k$.
\par Then we form a factorization $(K,g)$ of $f$ with arms $K_j=\al_j\cup D_j$, where disks $D_j$ are so small that $f$ is either holomorphic or antiholomorphic on them.  By Theorem \ref{T:rho} $[f_j]_\rho=[h_j]^{-1}_\rho$, where $h_j\in\A_{w_0}(K_j,W,M)$ when $n<j\le k$ and by Theorem \ref{T:ft}
\[[f]_\rho=\prod_{j=1}^n\dl_1([f_j])\prod_{j=n+1}^k(\dl_1([h_j]))^{-1}.\]
The part (4) has the same proof.
\end{proof}
\par By Theorem \ref{T:dls} if $\pi_1(M,w_0)=\pi_2(M,w_0)=0$ then the group $\rho_1$ in Theorem \ref{T:rr} can be replaced by $\pi_1(W,w_0)$.
\section{Connected components of $\A(W,M)$ and $\R(W,M)$}\label{S:cc}
\par There are natural mappings of the set $\eta_1(W,M)$ of connected components of $\A(W,M)$ into the set $\rho_1(W,M)$ of connected components of $\R(W,M)$ and  $\pi_1(W)$. We will denote these mapping also by $\dl_2$ and $\iota_1$ respectively. The mapping $\iota_1$ need not to be an injection especially when the group $\pi_2(M)$ is non-trivial. For example, if $M={\mathbb{CP}}^n$ and $A$ is an algebraic variety in $M$ such that $\pi_1(W)$ is finite then $\iota_1$ is not an injection because $\eta_1(W,M)$ is always infinite due to the invariance of index.
\par There is hope that $\dl_2$ is an injection at least when $M=\aC^n$. To advance in this direction we introduce the set $\R^\pm_{w_0}(W,M)$ of mappings $f\in \R_{w_0}(W,M)$ such that the set $A(f)$ is finite and the points in $A(f)$ have non-zero local indexes. This set is open and dense in $\R_{w_0}(W,M)$.
\bL\label{L:tl} Let $M$ be a complex manifold and let $A$ be an analytic variety in $M$ and $W=M\sm A$. Suppose that  mappings $f_0,f_1\in\R^\pm_{w_0}(W,M)$ are smooth and transverse to $A$ and can be connected by a continuous curve $f_t$, $0\le t\le1$, in $\R^\pm_{w_0}(W,M)$. Then there is  a smooth path $g_t$, $0\le t\le1$, connecting $f_0$ and $f_1$ in $\R^\pm_{w_0}(W,M)$ such that:\be
\item the mapping $G(t,\zeta)=g_t(\zeta)$ of $[0,1]\times\oD$ into $M$ is transverse to $A$;
\item the set $A_G$ consists of finitely many disjoint smooth curves  $\{\al_j(t)\}$, $t\in[0,1]$, such that $\\al_j(t)=(t,\zeta_j(t))$.
\ee \eL
\begin{proof} We may assume that the mapping $F:\,[0,1]\times\oD\to M$ defined as $F(t,\zeta)=f_t(\zeta)$ is smooth. Let $A_{\sing}=A^1_{\sing}$ be the set  of singular points of $A$.  Define by induction the sets $A^k_{\sing}=(A^{k-1}_{\sing})_{\sing}$. For some $k\le n+1$ the set $A^k_{\sing}$ is empty and therefore the set $A^{k-1}_{\sing}$ is a manifold. By the Thom Transversality Theorem we can approximate $F$ by a smooth mapping $F_k$ transverse to $A^{k-1}_{\sing}$. By the definition of transversality $F_k([0,1]\times\oD)$ never meets the set $A^{k-1}_{\sing}$ if $\dim A^{k-1}_{\sing}\le n-2$. Now we let $M_{k-1}=M\sm A^{k-1}_{\sing}$ and apply the transversality theorem to $M_{k-1}$ and $A^{k-2}_{\sing}$ to find $F_{k-2}$. By induction we obtain an approximation $H$ of $F$ that never meets the set $A_{\sing}$ and is transverse to $A$. Let $h_t(\zeta)=H(t,\zeta)$. Since $M$ admits a real analytic imbedding into some $\aR^N$ we can choose $H$ to be real analytic and since the set $\R^\pm_{w_0}(W,M)$ is open we may assume that $h_t\in\R^\pm_{w_0}(W,M)$ for all $t\in[0,1]$.
\par  The set $A_H$ is a compact set in $[0,1]\times\aD$ and a smooth submanifold, i.e., it is a collection $\Gm$ of finitely many disjoint smooth curves  $\{\gm_j\}$, $1\le j\le m$.
\par Suppose that $(t_0,\zeta_0)\in\gm_j$ and $H(t_0,\zeta_0)=w_0\in A$. The point $w_0$ is a regular point of $A$ and there is a neighborhood $U$ of $w_0$ such that in appropriate coordinates $(z_1,\dots,z_n)$ the set $A\cap U=\{z_1=0\}$. In coordinates $(z_1,\dots,z_n)$ the mapping $H(t,\zeta)=(H_1(t,\zeta),\dots,H_n(t,\zeta))$ and the functions $H_k$ are real analytic. Since the rank of $dH_1$ is 2 and the curve $\gm_j=\{H_1(t,\zeta)=0\}$, by the Implicit Function Theorem the curve $\gm_j$ admits a real analytic parametrization $\gm_j(s)=(t(s),\zeta(s))$ near $(t_0,\zeta_0)$ with $t(0)=t_0$. The mapping $h_{t_0}$ is not transverse to $A$ at $\zeta_0$ if and only if either $t(s)=t_0$ near $0$ or $t(s)=t_0+as^p+o(s^p)$, $a\ne0$ and $p>1$. But the former case is excluded because $h_{t_0}\in\R^\pm_{w_0}(W,M)$ and the set $A(h_{t_0})$ is finite. By real analyticity in the latter case there are only finitely many points where $t'(s)=0$. Hence the set $E$ of those $t$ where $h_t$ is not transverse to $A$ is finite.
\par Let $t_0\in E$. If $p$ is even and $a>0$ then $t(s)$ has a strict local minimum at $0$ and if $a<0$ then it has a strict local maximum there. In both cases $\ind_{\zeta_0}(h_{t_0},A)=0$ and this contradicts the assumption that $h_{t_0}\in\R^\pm_{w_0}(W,M)$. If $p>1$ is odd then $t(s)$ is either strictly increasing or decreasing near 0 and the set $h_t(\oD)\cap A$ has only one point in a small neighborhood of $\zeta_0$ for $t$ sufficiently close to $\zeta_0$.  Hence we can choose a real analytic parametrization $\gm_j(s)$ such that $s\in[0,1]$, the function $\tau_j(s)=t(\gm_j(s))$ is strictly increasing and $\tau_j'(s)=0$ only at finitely many points.
\par Finally, we take smooth diffeomorphism $\Phi(t,\zeta)$ of $[0,1]\times\oD$ such that $\Phi(t,\zeta)=(t,\phi(t,\zeta))$ and the functions $t(\Phi(\gm_j(s)))$ have strictly positive derivatives. Let $s_j(t)$ be the inverse of the later function. The mapping $G=H\circ\Phi^{-1}$ and the curves $\al_j(t)=\Phi(\gm_j(s_j(t))$ have all required properties.
\end{proof}
\par The proof of the lemma below follows the same line of argument as that in the proof of Assertion 2 in the proof of \cite[Lemma 2.1]{Sl}.
\bL\label{L:dif} Let $\zeta_k(t)$, $1\le k\le n$, are smooth mappings of $[0,1]$ into $\aD$ such that $\zeta_i(t)\ne\zeta_j(t)$ when $i\ne j$ and $0\le t\le 1$.  Then there is a $C^{\infty}$ mapping $\Phi: \ovr{\mathbb D}\times [0,1] \to \ovr{\mathbb D}$ such that $\Phi_t(\zeta)=\Phi(\zeta,t)$ is a diffeomorphism of  $\ovr {\mathbb D}$ onto itself for each $t$, $\Phi_t(\zeta)=\zeta$ when $|\zeta|=1$ and $\Phi_t(\zeta_j(0))=\zeta_j(t)$ for $j=1,\dots,n$.
\eL
\begin{proof} By Whitney extension theorem (see \cite[Theorem 1.5.6]{N}) we can find
a complex valued $C^{\infty}$-function $F(t,\zeta)$ on $[0,1]\times{\mathbb C}$ such that
$F(t,\zeta_j(t))=\partial\zeta_j(t)/\partial t$ for $0 \le t \le 1$, $j=1, \cdots, n$.
Replacing $F$ with the product $F\phi$, where $\phi$ is a $C^{\infty}$-function with $\phi=1$ on $\cup_{j=1}^n \{(t,\zeta_j(t)):\,0 \le t \le 1\}$ and $\phi=0$ for $|\zeta| \ge 1$, we can make $F(t,\zeta)=0$ for $|\zeta| \ge 1$. Then by standard existence and uniqueness theorems for ordinary differential equations, the initial value problem $\partial x/\partial t(t)=F(t,x(t))$, $x(0,\zeta)=\zeta$, $0 \le t \le 1$, has a unique solution $x(t,\zeta)$. Since $F(t,\zeta)$ is smooth, this solution  is smooth on $[0,1]\times{\mathbb C}$.
\par Now define a mapping $\Phi: {\mathbb C}\times [0,1] \to {\mathbb C}$ by $\Phi(\zeta,t)=x(t,\zeta)$. Then for each $0 \le t \le 1$, $\Phi_t$ is a diffemorphism and since the related initial value problem has a unique solution  we have $\Phi(\zeta_j(0),t)=\zeta_j(t)$ for $j=1,\cdots,n$. Also note that $\Phi(\zeta,t)=\zeta$ for all $0 \le t \le 1$ when $|\zeta| \ge1$. So, the restriction of $\Phi$ to $\ovr {\mathbb D} \times [0,1]$ has desired properties.
\end{proof}
\par Let $\gm:\,[0,1]\to M$ be a continuous curve connecting  the base point $w_0$ with a point $w$ in a regular part of some component $A_i$ and such that $\gm([0,1))\sbs W$. Given a neighborhood $U$ of $w$ and $\eps>0$ consider continuous mappings $f_{U,\eps}$ of $[0,1-\eps]\cup\oD(2-\eps,1)$ such that $f(t)=\gm(t)$ when $0\le t\le 1-\eps$ and the restriction of $f$ to $\oD(2-\eps,1)$ is an analytic disk in $U$ transversal to $A_i$ and whose index with respect to $A_i$ is 1. Rephrasing the definition from \cite{Sh} we call such mappings {\it lassos} $\lm_\gm$ around $A_i$. Clearly, there are $U$ and $\eps_0$ such that all mappings $f_{V,\eps}$ are equivalent in $\H[\oD,W,M]$ when $V\sbs U$ and $\eps<\eps_0$.
\bL\label{L:hl} Let $\lm_{\gm_0}$ and $\lm_{\gm_1}$ be lassos around $A_i$. Suppose that there is a continuous mapping $\phi:\,[0,1]^2\to M$ such that for all $t$ we have $\phi(0,t)=\gm_0(t)$, $\phi(1,t)=\gm_1(t)$, $\phi(t,0)=w_0$, $\phi(t,1)\in A^{\reg}_i$ and $\phi(t,s)\in W$ when $s\ne1$. Then $\lm_{\gm_0}$ and $\lm_{\gm_1}$ are equivalent in $\H_{0,w_0}[\oD,W,M]$
\eL
\begin{proof} For some small $\eps>0$ we can construct a continuous family $g_t$ of analytic disks transversal to $A_i$ and of index 1, centered at $\phi(t,1)$ and such that $g_t(1)=\phi(t,1-\eps)$. Then define $K=[0,1-\eps]\cup\oD(2-\eps,1)$ and $\lm_{\gm_t}:\,K\to M$ as $\phi(t,s)$ on $[0,1-\eps]$ and as $g_t(\zeta-2+\eps)$ on $\oD(2-\eps,1)$. The path $\lm_{\gm_t}$ is continuous in $\T(W,M)$ and by Theorem \ref{T:ci} $[\lm_{\gm_0}]_{0,w_0}=[\lm_{\gm_1}]_{0,w_0}$.
\end{proof}
\bT Let $M$ be a complex manifold and let $A$ be an analytic set in $M$ and let $W=M\sm A$. If $f_0$ and $f_1$ in $\A_{1,w_0}(W,M)$ belong to the same connected component of $\R^\pm_{w_0}(W,M)$ then $[f_0]_{1,w_0}=[f_1]_{1,w_0}$. \eT
\begin{proof} We may assume that $f_0$ and $f_1$ are transverse to $A$. Since they belong to the same connected component of $\R^\pm_{w_0}(W,M)$ there is a smooth mapping $G(t,\zeta):\,[0,1]\times\oD\to M$ satisfying conclusions of Lemma \ref{L:tl}. Let us apply Lemma \ref{L:dif} to the curves $\gm_j(t)$ in Lemma \ref{L:tl} to get the diffeomorphisms $\Phi_t$ of $\oD$.
\par Let us take a factorization $(K,g_0)$ of $f_0$ consisting of arms $(K_j,g_j)$, $1\le j\le k$, such that $f_0(K_j)$  are lassos $\al_j$ so that
\[ [f_0]=\prod_{j=1}^{k}[\al_j].\]
\par Let $L_j=\Phi_1(K_j)$. Then $L_j$ form a star and mappings $f_1(L_j)$ are also lassos. Moreover, the lassos $\al_j$ and $\beta_j$ are $h$-homotopic. The needed $h$-homotopy is achieved by the family $f_t(\Phi_t(K_j))$. Thus $[\al_j]=[\beta_j]$. Since
\[ [f_1]=\prod_{j=1}^{k}[\beta_j]\] the theorem is proved.
\end{proof}

\end{document}